 \documentclass[francais,brochure]{bourbaki}

 \usepackage[T1]{fontenc}
 \usepackage{xspace,amssymb,mathrsfs,url}
 
 \usepackage[frenchb]{babel}
 \usepackage{ifpdf}
 \ifpdf
      \usepackage{ae}\providecommand{\og}{``}\providecommand{\fg}{''}
 \fi

\date{Mars 2003}
\bbkannee{55\textsuperscript{e} ann{\'e}e, 2002-2003}
\bbknumero{914}

\title
 [Quelques applications des cohomologies $p$-adiques]
 {Points rationnels et groupes fondamentaux : \\
  applications de la cohomologie $p$-adique}
\subtitle{d'apr\`es P.~Berthelot, T.~Ekedahl, H.~Esnault, etc.}
\author{Antoine Chambert-Loir}
\address{Centre de math\'ematiques (CMAT) \\ \'Ecole polytechnique \\ 
91128 Palaiseau Cedex}
\email{chambert@math.polytechnique.fr}

\def\theenumi{\alph{enumi}}\def\labelenumi{\theenumi)}
\makeatletter
  \let\c@equation\c@defi
  \let\theequation\thedefi
  \let\cl@equation\cl@defi
\makeatother

\def\card#1{\mathopen| #1\mathclose|}
\def\A{\mathbf{A}}
\def\P{\mathbf{P}}
\def\F{\mathbf{F}}
\def\Z{\mathbf{Z}}
\def\C{\mathbf{C}}
\def\R{\mathbf{R}}
\def\Q{\mathbf{Q}}
\def\Nwt{\operatorname{Nwt}}
\def\Tr{\operatorname{Tr}}
\def\Hdg{\operatorname{Hdg}}
\def\tube#1{\mathopen]#1\mathclose[}

\let\phi\varphi

\def\ieme{\nobreakdash-i\`eme\xspace}
\def\iemes{\nobreakdash-i\`emes\xspace}
\def\cf{\emph{cf.}~}
\def\resp{\emph{resp.}~}
\let\ra\rightarrow
\def\sozat{\;;\;}
\def\CH{\operatorname{CH}}
\def\Gal{\operatorname{Gal}}

\def\Frac{\operatorname{Frac}}
\def\alb{\operatorname{alb}}
\def\End{\operatorname{End}}
\def\id{\operatorname{id}}
\def\ch{\operatorname{ch}}
\def\norm#1{\mathopen\|{#1}\mathclose\|}
\def\abs#1{\mathopen|{#1}\mathclose|}
\def\ord{\operatorname{ord}}
\def\rig{{\text{\upshape rig}}}
\def\etc{{\et,\mathrm c}}
\def\et{{\text{\upshape \'et}}}
\def\cris{{\text{\upshape cris}}}
\def\rigc{{\rig,\mathrm c}}

\begin{document}

\maketitle

Dans cet expos\'e, je pr\'esente 
trois r\'esultats concernant les vari\'et\'es alg\'ebriques
en caract\'eristique positive :

\begin{enumerate}\itshape
\item 
Deux vari\'et\'es propres et lisses sur $\F_q$ 
qui sont g\'eom\'etriquement birationnelles
ont m\^eme nombre de points rationnels modulo~$q$ 
\emph{(\cf T.~Ekedahl, \cite{ekedahl1983})}.

\item
Sur un corps fini,
une vari\'et\'e propre et lisse qui est de Fano, 
ou bien g\'eom\'etriquement faiblement unirationnelle, 
ou plus g\'en\'eralement rationnellement 
connexe par cha\^{\i}nes, a un point rationnel 
\emph{(H.~Esnault, \cite{esnault2003})}.

\item
Sur un corps alg\'ebriquement clos de caract\'eristique~$p>0$,
le groupe fondamental d'une vari\'et\'e propre et lisse
qui est de Fano, ou bien g\'eom\'etriquement
faiblement unirationnelle, ou plus g\'en\'eralement rationnellement 
connexe par cha\^{\i}nes, est un groupe fini
d'ordre premier \`a~$p$
\emph{(\cf T.~Ekedahl, \cite{ekedahl1983})}.
\end{enumerate}

Le point commun des d\'emonstrations est un contr\^ole des valuations
$p$-adiques des valeurs propres de Frobenius. Elles gagnent donc
\`a \^etre pr\'esent\'ees dans le cadre d'une th\'eorie cohomologique $p$-adique.
La cohomologie rigide, d\'evelopp\'ee par P.~Berthelot, fournit l'outil
id\'eal pour cela. Elle a connu r\'ecemment des progr\`es importants 
et je d\'ecris le formalisme auquel elle donne lieu.
Les \'enonc\'es ci-dessus s'obtiennent en contr\^olant
les \emph{pentes} des F-isocristaux que fournit la cohomologie rigide.

\bigskip
{\itshape
Je voudrais remercier 
P.~Berthelot, J.-L.~Colliot-Th\'el\`ene, O.~Debarre, H.~Esnault, B.~Kahn,
Y.~Laszlo et J-P.~Serre pour leurs conseils ou suggestions.
\par}

% \part{Amuse-bouches}

\section{Autour du th\'eor\`eme de Chevalley-Warning}

Je commence cet expos\'e par un \'enonc\'e \'el\'ementaire
et assez ancien, d\^u \`a C.~Chevalley et E.~Warning~\cite{warning36}.

\begin{theo}\label{theo.cw}
Soit $\F$ un corps fini, $q$ son cardinal et $p$ sa caract\'eristique.
Soit $F_1,\dots,F_r$ des polyn\^omes en $x_1,\dots,x_n$,
\`a coefficients dans~$\F$ et de degr\'es $d_1,\dots,d_r$.
Si $n>d_1+\dots+d_r$,
le nombre de solutions dans $\F^n$ du syst\`eme
\begin{equation}
\label{eq.cw}
 F_1(x_1,\dots,x_n) = \dots=F_r(x_1,\dots,x_n) =0
\end{equation}
est multiple de~$p$.
\end{theo}

La d\'emonstration classique est tr\`es simple et repose sur le fait
que pour $x\in\F$, $x^{q-1}$ vaut $1$ si $x\neq 0$ et $0$ sinon.
Ainsi, le nombre de solutions du syst\`eme est congru modulo~$p$
\`a l'expression
\[ \sum_{\mathbf x\in\F^n} \prod_{i=1}^r (1-F_i(\mathbf x)^{q-1}) .\]
Le produit 
\[ \prod_{i=1}^r (1-F_i(\mathbf x)^{q-1}) \]
est un polyn\^ome de degr\'e $\leq (q-1)\sum d_i$. Soit 
\[ c_{\mathbf m} x_1^{m_1}\dots x_n^{m_n} \]
un de ses mon\^omes non nuls. On a donc $m_1 +\dots + m_n\leq
(q-1)\sum d_i < n(q-1)$
si bien que n\'ecessairement, l'un des entiers $m_i$, disons $m_1$,
est strictement
inf\'erieur \`a $q-1$ et 
\[  \sum_{\mathbf x\in\F^n} c_{\mathbf m} x_1^{m_1}\dots x_n^{m_n}
  = c_{\mathbf m} \prod_{i=1}^n \sum_{t\in\F} t^{m_i}
 = 0 \]
puisque
\[ \sum_{t\in\F} t^{m} = \begin{cases} -1 & \text{si $(q-1)$ divise $m$ et
$m\geq 1$,} \\
      0 &\text{sinon.} \end{cases}
\]

Ce th\'eor\`eme a \'et\'e g\'en\'eralis\'e par J.~Ax~\cite{ax1964} et N.~Katz~\cite{katz1971}:
\begin{theo}\label{theo.ak}
Si $b$ d\'esigne le plus petit entier tel que
\[ b \geq \frac{n-\sum d_i}{\max d_i}, \]
le nombre de solutions du syst\`eme~\eqref{eq.cw}
est divisible par $q^b$.
\end{theo}

Leurs d\'emonstrations sont assez savantes. 
Celle qu'a propos\'ee r\'ecemment D.~Wan dans~\cite{wan1995} est en revanche
\'el\'ementaire et tout \`a fait dans l'esprit de celle
du th\'eor\`eme~\ref{theo.cw}. Elle commence par une r\'eduction des
scalaires (\`a la Weil) au cas du corps premier $\F_p$,
introduite dans ce contexte par C.~Moreno et O.~Moreno~\cite{moreno-m1995}:
si on fixe une base de $\F$ sur $\F_p$, de cardinal $a$, on se ram\`ene
\`a un syst\`eme de $ar$ \'equations en $nr$ variables de degr\'es
$d_1,\dots,d_r$ ($a$ fois). Le nouvel entier $b$ est ainsi \'egal \`a
$a$ fois l'ancien si bien que la congruence pour le nouveau syst\`eme
implique celle du th\'eor\`eme~\ref{theo.ak}.

Pour continuer la d\'emonstration, dans le cas o\`u $\F=\F_p$,
il est commode d'introduire l'anneau $\Z_p$ des entiers $p$-adiques.
C'est un anneau de valuation discr\`ete complet, de caract\'eristique~$0$,
de corps r\'esiduel $\F_p$ ; son id\'eal maximal est engendr\'e par $p$ ;
il contient les $p$ racines de l'\'equation $x^p-x=0$ (notons $S$ cet
ensemble).
L'id\'ee de base est maintenant que pour $x\in \Z_p$, 
$ x^{p^n(p-1)}$ est congru modulo $p^n$ \`a $1$ si $x\not\equiv 0 \pmod p$,
et \`a $0$ si $x\equiv 0\pmod p$.
Le nombre de solutions du syst\`eme~\eqref{eq.cw} est ainsi congru modulo~$p^n$
\`a l'expression
\[ \sum_{\mathbf x\in S^n} \prod_{i=1}^r \big(1-F_i(\mathbf x)^{p^n(p-1)}\big)
.\]
Par r\'ecurrence sur $r$, il suffit de montrer que
\[ \sum_{\mathbf x\in S^n} \prod_{i=1}^r F_i(\mathbf x)^{p^n(p-1)}
 \equiv 0  \pmod{p^b}, \]
ce que fait Wan en d\'eveloppant $F_i^{p^n(p-1)}$ puis en constatant
que la valuation des divers symboles du multin\^omes et celle de la somme
de puissances qui apparaissent s'ajoutent pour d\'epasser $b$.
(C'est tout de m\^eme assez d\'elicat.)

\section{Fonctions z\^eta et cohomologies de Weil}

Soit $\F$ un corps fini, notons~$q$ son cardinal.
Dans~\cite{weil49}, A.~Weil avait introduit pour tout syst\`eme d'\'equations
polyn\^omiales \`a coefficients dans~$\F$, voire
tout $\F$-sch\'ema $V$ de type fini, la s\'erie
g\'en\'eratrice 
\[ Z_V(t) = \exp \left( \sum_{a\geq 0} \card{ V(\F^{(a)})} \frac{t^a}a \right)
\]
o\`u $\F^{(a)}$ d\'esigne l'unique extension de $\F$ de degr\'e~$a$.
B.~Dwork~\cite{dwork60} a d\'emontr\'e que c'est une fraction rationnelle 
(premi\`ere conjecture de Weil).
Ses z\'eros et ses p\^oles sont n\'ecessairement des entiers alg\'ebriques.
En fait, 
la congruence du th\'eor\`eme~\ref{theo.ak} implique
qu'ils sont divisibles par $q^b$ dans l'anneau des entiers
alg\'ebriques.

Les deux autres conjectures de Weil ont \'et\'e d\'emontr\'ees
par A.~Grothendieck (rationalit\'e et \'equation fonctionnelle)
dans SGA~5
et P.~Deligne (analogue de l'hypoth\`ese de Riemann, \cite{deligne74}),
gr\^ace \`a l'introduction de la cohomologie (\'etale) $\ell$-adique,
o\`u $\ell$ est un nombre premier \emph{distinct} de $p$.
L'utilit\'e d'une telle th\'eorie cohomologique avait d\'ej\`a \'et\'e pressentie
dans l'article de Weil. 
%% RMK : pas gentil ! (JPS)
%% --- pour lequel il avait command\'e \`a J.~Dolbeault
%% le calcul des nombres de Betti  des hypersurfaces de l'espace projectif.
En effet, $V$ dispose d'un endomorphisme de Frobenius $F$
et l'ensemble $V(\F^{(a)})$ n'est autre  que le lieu des points
fixes de $F^a$ ; il faudrait pouvoir disposer alors d'une
formule des traces de Lefschetz.
Les conditions n\'ecessaires  sur une telle cohomologie
ont \'et\'e rapidement formalis\'ees sous le nom
de \emph{cohomologie de Weil} : la th\'eorie cohomologique
doit v\'erifier la formule de K\"unneth,
la dualit\'e de Poincar\'e, fournie par une classe fondamentale,
et les sous-vari\'et\'es (non n\'ecessairement lisses) doivent
poss\'eder une classe de cohomologie, compatibles au produit
d'intersection et \`a la classe fondamentale.
Tout ceci est expos\'e en d\'etail dans~\cite{kleiman68}, 
en lien avec les th\'eor\`emes
de Lefschetz faible et difficile, 
les conjectures {\og standard\fg}
et la conjecture de Hodge.

Lorsque $\ell\neq p$, la cohomologie \'etale $\ell$-adique
est effectivement une cohomologie de Weil.
Lorsque $\ell=p$, ce n'est pas vrai : $H^i(X,\Q_p)$ est
nul d\`es que $i>\dim X$ ce qui viole la dualit\'e de Poincar\'e.

Finalement, par voie $\ell$-adique, on sait que la fonction
z\^eta d'un $\F$-sch\'ema
g\'eom\'etriquement connexe, propre, lisse de dimension~$d$, 
est de la forme suivante :
\[ Z_X (t) = \frac{P_1(t)\dots P_{2d-1}(t)}{P_0(t)\dots P_{2d}(t)}, \]
o\`u $P_i(t)\in\Z[t]$, avec $P_0(t)=1-t$, $P_{2d}(t)=1-q^{2d}t$ ;
si $P_{i}(t)=\prod_{j=1}^{b_i} (1-a_{ij}t)$, on sait aussi
que les $a_{ij}$ sont des entiers alg\'ebriques de valeur
absolue archim\'edienne $q^{i/2}$ (hypoth\`ese de Riemann)
et, si $\ell\neq p$, de valeur absolue $\ell$-adique nulle.

L'interpr\'etation des valuations $p$-adiques des $a_{ij}$,
et notamment de congruences du type fourni par les th\'eor\`emes~\ref{theo.cw}
et~\ref{theo.ak} justifie la recherche d'une th\'eorie cohomologique
de Weil qui soit $p$-adique.

Il y a un autre int\'er\^et 
--- sur lequel insiste Kedlaya ---
des cohomologies $p$-adiques que je vais d\'ecrire
maintenant :
elles sont par nature plus calculables
que ne l'est la cohomologie \'etale. Il est par exemple remarquable
que soient apparus r\'ecemment trois \emph{algorithmes} efficaces
(Satoh~\cite{satoh2000},  Kedlaya~\cite{kedlaya2001},
 Lauder-Wan~\cite{lauder-wan2002})
pour calculer le nombre de points de certaines vari\'et\'es
alg\'ebriques sur un corps fini, et tous trois sont de nature $p$-adique.
Celui de Kedlaya repose sur l'explicitation
de la cohomologie de Monsky-Washnitzer, celui de
Lauder-Wan
est inspir\'e de la d\'emonstration par Dwork de la rationalit\'e de la fonction
z\^eta.
De m\^eme, c'est par des techniques $p$-adiques que
Bombieri~\cite{bombieri1978} a \'etudi\'e le degr\'e des num\'erateurs et
d\'enominateurs de fonctions rationnelles associ\'ees \`a des sommes
d'exponentielles.

\section{Cohomologies $p$-adiques}

Mais avant cela, il faut peut-\^etre dire pourquoi il n'existe pas
de cohomologie de Weil sur la cat\'egorie
des vari\'et\'es alg\'ebriques projectives lisses
sur un corps alg\'ebriquement clos de caract\'eristique
positive \`a valeur dans la cat\'egorie des $\Q$-espaces vectoriels 
de dimension finie. Pour toute cohomologie de Weil \`a valeurs
dans les $K$-vectoriels ($K$ est un corps commutatif),
le $H^1$ d'une courbe elliptique~$E$ est de dimension~$2$.
De plus, pour tout endomorphisme non nul $\alpha$  de~$E$,
l'endomorphisme $\alpha^*\colon H^1(E)\ra H^1(E)$ est injectif,
d'o\`u une injection de l'anneau (oppos\'e \`a celui) des endomorphismes de $E$
dans l'anneau des matrices $2\times 2$ \`a coefficients dans~$K$.
En caract\'eristique~$0$, tout irait bien (d'ailleurs, la cohomologie
singuli\`ere est une cohomologie de Weil),
mais en caract\'eristique finie, il existe des courbes elliptiques
\emph{supersinguli\`eres} dont l'anneau des endomorphismes 
est plus gros que pr\'evu : c'est une alg\`ebre de quaternions sur~$\Q$.
Comme une telle alg\`ebre ne se plonge pas dans $M_2(\Q)$, il n'y
a pas de cohomologie de Weil \`a coefficients dans~$\Q$.
M.~Deuring a m\^eme montr\'e que 
l'alg\`ebre de quaternions $\End(E)\otimes\Q$ est ramifi\'ee en $p$ et l'infini,
c'est-\`a-dire que
$\End(E)\otimes\Q_p$ et $\End(E)\otimes\R$ sont des corps gauches.
Cela emp\^eche aussi $K=\Q_p$ ou $K=\R$.
% RMK (JPS)
% (Cet argument est tir\'e de~\cite{grothendieck68}, o\`u il est attribu\'e \`a 
% J-P.~Serre.)
(Cet argument est d\^u \`a J-P.~Serre, \cf\cite{grothendieck68}.)

La th\'eorie que nous utiliserons dans cet expos\'e est la \emph{cohomologie
rigide}, construite par P.~Berthelot. Elle unifie deux th\'eories disjointes
qui sont la cohomologie de Monsky-Washnitzer~\cite{monsky-w68},
valable pour les vari\'et\'es affines et lisses, et la cohomologie
cristalline~\cite{grothendieck68,berthelot74b}, qui a de bonnes
propri\'et\'es pour les vari\'et\'es propres et lisses.
Cette th\'eorie est encore \'eparpill\'ee dans la litt\'erature, 
et un petit guide de lecture ne sera peut-\^etre pas inutile.
Il faudrait aussi citer les travaux de Y.~Andr\'e, 
F.~Baldassarri, P.~Berthelot, B.~Chiarellotto, G.~Christol, 
R.~Crew, B.~Dwork, Z.~Mebkhout, P.~Robba, N.~Tsuzuki...

Je trouve l'introduction que propose Berthelot dans \cite{berthelot1988}
tr\`es agr\'eable \`a lire ; le gros article~\cite{berthelot1996}
fournit des d\'etails importants dans la construction 
(les fameux {\og th\'eor\`emes de fibrations\fg}).
L'article~\cite{berthelot1997} 
est tr\`es important : outre la d\'emonstration du fait que la cohomologie
rigide d'une vari\'et\'e lisse est de dimension finie, 
on y trouve les th\'eor\`emes de comparaison
avec les th\'eories de Monsky-Washnitzer et cristalline.
Concernant ce th\'eor\`eme de finitude, citons aussi l'article~\cite{mebkhout1997}
de Z.~Mebkhout
qui propose une d\'emonstration de la finitude de la cohomologie
de Monsky-Washnitzer ind\'ependante des th\'eor\`emes
de de Jong sur l'existence d'alt\'erations. Notons
que la finitude de la cohomologie de Monsky-Washnitzer n'\'etait pas connue
avant ces deux articles, \`a l'exception du $H^0$ et du $H^1$ par
P.~Monsky.
La finitude de la cohomologie rigide dans le cas g\'en\'eral est
d\'emontr\'ee (ind\'ependamment) 
dans l'article~\cite{grosse-klonne2002} de E.~Grosse-Kl\"onne
et dans celui~\cite{tsuzuki2003} de N.~Tsuzuki dans lequel
il \'etablit un th\'eor\`eme de descente cohomologique propre.

La dualit\'e de Poincar\'e et la formule de K\"unneth
sont \'etablies dans la note~\cite{berthelot1997b} de Berthelot.
Enfin, l'existence de classes de Chern est d\'emontr\'ee par
D.~Petrequin~\cite{petrequin2003}.

\bigskip

\emph{Qu'est ce que la cohomologie rigide ?}
Disons tout d'abord que c'est une sorte de cohomologie de de Rham.
Fixons quelques notations : soit $k$ un corps 
de caract\'eristique $p>0$, fixons alors un anneau de valuation discr\`ete 
complet $\mathcal V$ de corps r\'esiduel $k$,
dont nous noterons $K$ le corps des fractions, suppos\'e
de caract\'eristique z\'ero, et $\pi$ un g\'en\'erateur
de l'id\'eal maximal de $\mathcal V$. 
Nous supposerons que $\mathcal V$
admet un endomorphisme $\sigma$ qui se r\'eduit modulo~$\pi$
en l'automorphisme de Frobenius $\sigma\colon x\mapsto x^p$ de $k$ ;
alors, $\sigma$ s'\'etend en un endomorphisme de $K$. Dans tout
ce qui va suivre, on peut se limiter au cas important
o\`u le corps~$k$ est suppos\'e parfait et o\`u $\mathcal V$ est l'anneau 
des vecteurs de Witt de~$k$.

Soit $X$ un $k$-sch\'ema et essayons de d\'efinir la cohomologie
rigide de $k$. Ce seront des $K$-espaces vectoriels.
Le cas id\'eal est celui o\`u $X$ est la \emph{r\'eduction} modulo~$p$
d'un $V$-sch\'ema propre et lisse $\mathcal X$. Dans ce cas, 
$\mathcal X$ dispose d'une cohomologie de de Rham d\'efinie alg\'ebriquement :
l'hypercohomologie du complexe des formes diff\'erentielles sur $\tilde X$;
ce sont des $\mathcal V$-modules de type fini.
Un point crucial, d\'ej\`a \`a la base de l'existence de la cohomologie
de Monsky--Washnitzer, est que ces modules ne d\'ependent pas du choix
de $\mathcal X$ --- autrement dit, si $\mathcal X'$ est un autre
$\mathcal V$-sch\'ema propre et lisse de m\^eme r\'eduction modulo~$p$,
on a un isomorphisme canonique 
$H^*_{DR}(\mathcal X)\simeq H^*_{DR}(\mathcal X')$.

Du point de vue topologique, une vari\'et\'e plong\'ee dans un espace
lisse est r\'etracte par d\'eformation d'un voisinage tubulaire assez petit,
et il est possible de d\'efinir sa cohomologie \`a l'aide de celle
de ces voisinages.
En caract\'eristique z\'ero, Hartshorne~\cite{hartshorne1975} avait
\'etudi\'e une cohomologie de de Rham pour des vari\'et\'es singuli\`eres
d\'efinie suivant ces lignes ;  le r\^ole du voisinage tubulaire y est jou\'e
par le compl\'et\'e formel de l'espace ambiant le long de la vari\'et\'e 
dont il s'agit de d\'efinir la cohomologie.

La d\'efinition de la cohomologie rigide {\og na\"{\i}ve\fg} suit
cette approche si ce n'est qu'il faut plonger et relever.
Supposons donc que $X$ est un sous-sch\'ema d'un sch\'ema propre et lisse $P$,
r\'eduction d'un $\mathcal V$-sch\'ema $\mathcal P$ ; l'exemple important
est bien entendu l'espace projectif. \`A $\mathcal P$ est associ\'ee
une vari\'et\'e analytique \emph{rigide}, not\'ee $\mathsf P$ --- c'est une
structure plus riche que la structure analytique $p$-adique na\"{\i}ve
sur $\mathcal P\otimes K$ qui 
donne lieu \`a une th\'eorie des faisceaux non triviale
(penser que la topologie de $K$ est totalement discontinue !).

Un point de $\mathsf P$ se \emph{sp\'ecialise} en un point de $P$ :
dans le cas de l'espace projectif,
il suffit de chasser les d\'enominateurs pour qu'un point soit \`a
coordonn\'ees homog\`enes dans $\mathcal V$, non toutes multiples de~$\pi$,
puis de r\'eduire modulo~$\pi$ ces coordonn\'ees homog\`enes. Dans $\mathsf P$,
on peut alors d\'efinir le \emph{tube} de $X$ comme l'ensemble des points
de $\mathsf P$ qui se r\'eduisent en un point de~$X$. Ce tube, not\'e usuellement
$\tube X$, est une vari\'et\'e analytique rigide (pas forc\'ement
quasi-compacte).
Par exemple, si $X=\{0\}$ et $\mathcal P=\P^1$, $\tube X$ s'identifie
au {\og disque unit\'e ouvert\fg}, form\'e des $x\in K$
tels que $\abs x<1$. Si $X=\A^1$ et $\mathcal P=\P^1$, $\tube X$
est alors le {\og disque unit\'e ferm\'e\fg}, form\'e des $x\in K$
tels que $\abs x\leq 1$.

La cohomologie rigide na\"{\i}ve de $X$ est la cohomologie de de Rham
du tube de~$X$. D'apr\`es un th\'eor\`eme de fibration, elle ne d\'epend
pas du choix de~$\mathcal P$. Si l'on peut prendre $P=X$, c'est-\`a-dire
si $X$ est propre, lisse et relevable, la cohomologie d\'efinie
n'est autre que la cohomologie de de Rham tensoris\'ee par~$K$.
Si $X$ est propre et lisse, on retrouve la cohomologie cristalline de~$X$
tensoris\'ee par~$K$.
Ce sont en particulier des $K$-espaces vectoriels de dimension finie.

En revanche, si $X$ n'est pas propre, cela ne suffit pas.
Prenons l'exemple de la droite affine $X=\A^1$ 
et de son tube $\tube X=\{x\sozat \abs x\leq 1\}$. 
Le complexe de de Rham dont la cohomologie rigide na\"{\i}ve est la cohomologie
est donn\'e par
\[ K\{x\} \ra K\{x\}, \quad f=\sum a_n x^n\mapsto f'=\sum na_n x^{n-1} , \]
o\`u $K\{x\}$ est l'anneau des s\'eries enti\`eres \`a coefficients dans~$K$
dont les coefficients tendent vers~$0$ (de sorte qu'elles convergent sur le
disque ferm\'e).
On a bien $H^0_{\text{na\"{\i}f}}(\A^1)=K$, mais $H^1_{\text{na\"{\i}f}}(\A^1)$
n'est pas nul puisque
la s\'erie $f=\sum p^n x^{p^n-1}$ n'a pas de primitive dans $K\{x\}$.
En fait, $H^1_{\text{na\"{\i}f}}(\A^1)$ est m\^eme de dimension infinie.

Monsky et Washnitzer ont remarqu\'e que la situation s'arrange notablement
si l'on remplace l'anneau $K\{x\}$ par celui des fonctions qui convergent
dans un disque un peu plus gros que le disque unit\'e (on dit qu'elles
\emph{surconvergent}).
Notons $K\langle x\rangle^\dagger$ cet anneau : il est form\'e
des s\'eries $\sum a_n x^n$ telles que $\limsup \log\abs{a_n}/\log n<0$.
Le complexe de de Rham surconvergent de la droite affine a alors
la cohomologie attendue car si $f=\sum a_n x^n$ converge sur le disque
$\abs x\leq\lambda$, avec $\lambda>1$,
ses primitives convergent sur le disque ouvert
$\abs x<\lambda$, donc sur tout disque ferm\'e $\abs x\leq\lambda'$
avec $1<\lambda'<\lambda$.

C'est ainsi que pour d\'efinir la cohomologie rigide en g\'en\'eral,
il faut introduire ce que Berthelot
appelle des \emph{voisinages stricts} du tube (analogues
des disques ferm\'es $\abs x\leq\lambda$ pour le disque unit\'e)
et la cohomologie du complexe de de Rham form\'e des
formes diff\'erentielles surconvergentes, c'est-\`a-dire que
convergent dans un voisinage strict non pr\'ecis\'e de $\tube X$.

Outre la cohomologie rigide $H^i_\rig(X/K)$,
Berthelot d\'efinit aussi  une cohomologie \`a support $H^*_{\rig,Z}(X/K)$
(pour $Z\subset X$) et une cohomologie \`a support propre
$H^*_\rigc(X/K)$.
Ils sont de dimension finie (voir les r\'ef\'erences plus haut).
Ce sont les analogues alg\'ebriques de la cohomologie d'une paire
et de la cohomologie \`a support propre.

\medskip

Ces espaces de cohomologie sont compatibles \`a l'extension des scalaires 
si $K'$ est une extension isom\'etrique de $K$, 
d'anneau de valuation $\mathcal V'$, de corps r\'esiduel $k'$, 
il existe un isomorphisme canonique
\[ K'\otimes_K H^*_\rig(X/K) \xrightarrow\sim H^*_\rig(X'/K'), \]
o\`u $X'=k'\otimes_k X$.
(Voir~\cite{berthelot1997}, prop.~1.8, pour le cas d'une
extension finie, le cas g\'en\'eral est plus difficile
et est affirm\'e \`a la fin de~\cite{berthelot1997b}.) 

\medskip

Disons un mot des fonctorialit\'es dont disposent ces cohomologies.
La cohomologie rigide est naturellement contravariante pour
les morphismes de~$k$-sch\'emas.
La cohomologie \`a support propre n'est contravariante que
pour les morphismes propres, et est covariante pour les immersions
ouvertes.
Sur la cohomologie rigide des vari\'et\'es lisses, 
la dualit\'e de Poincar\'e permet d'en d\'eduire
une fonctorialit\'e covariante (avec un d\'ecalage de deux fois
la dimension) pour les morphismes
propres ({\og morphismes de Gysin\fg}).

Le morphisme de Frobenius $F_X$ n'est pas un morphisme
de $k$-sch\'ema, sauf si $k=\F_p$, mais il se factorise
en 
\[ X  \xrightarrow {F_{X/k}} k\otimes_{\sigma} X \ra  X \]
o\`u le premier morphisme est $k$-lin\'eaire
et le second est le morphisme de changement de base
par le Frobenius $\sigma$ de~$k$.
Gr\^ace \`a la compatibilit\'e des scalaires, 
si $\sigma$ est un endomorphisme de $\mathcal V$
qui induit le Frobenius modulo~$\pi$, on en d\'eduit un endomorphisme
$\sigma$-lin\'eaire de la cohomologie rigide :
\[ F\colon H^*_\rig(X/K) \ra H^*_\rig(X/K), \quad F(ax)=\sigma(a)F(x) , \]
et de m\^eme pour les cohomologies \`a support et \`a support propre.

Les espaces de cohomologie rigide
s'ins\`erent dans des suites exactes d'excision famili\`eres :
si $U$ est un ouvert de $X$ et $Z$ le ferm\'e compl\'ementaire,
on a des suites exactes
\[  H^i_\rigc(Z/K) \ra H^i_\rigc(X) \ra H^i_\rigc(U) \xrightarrow{[+1]} \]
et
\[ H^i_{\rig,Z}(X/K) \ra H^i_\rig(X) \ra H^i_\rig(U) \xrightarrow{[+1]} .\]
Celles-ci sont d'ailleurs \'equivariantes pour les divers
morphismes de Frobenius (voir~\cite{chiarellotto1998}, th.~2.4).

\medskip

Dans tout ceci, je n'ai en fait parl\'e que des {\og coefficients
constants\fg}. L'analogue des faisceaux localement constants
est fourni par les $F$-cristaux surconvergents : ce sont des fibr\'es vectoriels
sur un voisinage strict du tube munis d'une connexion int\'egrable
et d'une structure de Frobenius. 
La cohomologie \`a support propre d'un F-isocristal surconvergent
et la formule des traces de type Lefschetz
sont \'etudi\'ees dans les articles~\cite{etesse-ls1993,etesse-ls1997}
d'\'Etesse et Le Stum.
La finitude de la cohomologie rigide d'un F-isocristal surconvergent,
la dualit\'e de Poincar\'e et la formule de K\"unneth sont d\'emontr\'ees
dans la pr\'epublication~\cite{kedlaya2002} de K.~Kedlaya.

La th\'eorie des $\mathcal D$-modules
arithm\'etiques de Berthelot est cens\'ee fournir une cat\'egorie
de coefficients stable par les six op\'erations de Grothendieck,
mais \`a ma connaissance, ceci n'est pas encore d\'emontr\'e.

% \part{Plat de r\'esistance}

\section{F-isocristaux}

(Pour ce paragraphe, l'article~\cite{katz79b} de Katz est un \emph{must}.)
Supposons pour simplifier que
$k$ est un corps parfait de caract\'eristique~$p>0$,
$\mathcal V$ l'anneau des vecteurs de Witt de~$k$
et $K$ son corps des fractions.
Soit $\sigma$ l'automorphisme de $\mathcal V$ qui rel\`eve 
l'automorphisme de Frobenius de~$k$ ; il s'\'etend \`a~$K$.

\begin{defi}
Un F-isocristal sur~$K$ est un $K$-espace vectoriel de dimension finie
muni d'un endomorphisme $\sigma$-lin\'eaire injectif.
\end{defi}

Par exemple, soit $\alpha\in\Q^*$ un rationnel non nul, notons
$\alpha=r/d$ avec $(r,d)=1$ et $d>0$, et soit $M_\alpha=K^d$,
de base $(e_1,\dots,e_d)$,
muni de l'endomorphisme $\sigma$-lin\'eaire $F$ donn\'e par 
\[ F(e_1)=e_2,\dots,\quad F(e_{d-1})=e_d, \quad F(e_d)=p^r e_1. \]
De m\^eme, les espaces de cohomologie rigide d'un sch\'ema de type
fini sont naturellement des $F$-isocristaux (l'injectivit\'e n'est
pas \'evidente et sera \'etablie au paragraphe suivant.)

Si $(M,F)$ est un F-isocristal,
on peut exprimer $F$ dans une $K$-base de~$M$ \`a l'aide d'une matrice
$A\in M_n(K)$  ($n=\dim_K M$).
Il y a aussi une notion \'evidente de somme directe, de produit tensoriel,
ext\'erieur, sym\'etrique, d'homomorphisme de $F$-isocristaux. 

Si le corps $k$ est alg\'ebriquement clos,
la cat\'egorie des $F$-isocristaux a \'et\'e \'elucid\'ee par J.~Dieudonn\'e et
Yu.~Manin~\cite{manin63} : tout F-isocristal est somme directe
de $F$-isocristaux (simples) du type $M_\alpha$.
Cela permet de d\'efinir les \emph{pentes} d'un F-isocristal~$M$ :
ce sont les rationnels $\alpha$ tels que $M_\alpha$ soit un sous-objet
de~$M$. Si $M\simeq \bigoplus M_\alpha^{n_\alpha}$, la mutiplicit\'e
de la pente~$\alpha$ dans~$M$ est par d\'efinition \'egale \`a 
$n_\alpha \dim M_\alpha$. 

Si le corps~$k$ n'est pas alg\'ebriquement clos, les pentes
d'un F-isocristal sont celles du F-isocristal obtenu
apr\`es tensorisation par le corps des fractions de $W(\bar{\F_p})$.

Quel que soit le corps~$k$, pour tout rationnel $\alpha$,
on peut d\'efinir facilement le plus grand
sous-F-isocristal $M^{\geq\alpha}$ de $M$ dont les pentes sont $\geq \alpha$.
Fixons une base de~$M$, d'o\`u on d\'eduit une norme ultram\'etrique
$\norm{\cdot}$ sur~$M$.
L'ensemble $M^{\geq\alpha}$ des $x\in M$ tels que la suite
$(\norm{F^n(x)}p^{\alpha n})_n$ soit born\'ee
est un sous-$K$-espace vectoriel de~$M$, stable par~$F$ ; il ne
d\'epend pas de la base choisie.
Cela d\'efinit une filtration d\'ecroissante de~$M$ par des
sous-$F$-isocristaux, exhaustive (si les
coefficients d'une matrice de~$F$ sont de valuation $\geq r$,
$M=M^{\geq r}$, plus un raisonnement analogue pour $F^{-1}$).

Pour synth\'etiser les pentes d'un F-isocristal, il est commode
de d\'efinir son \emph{polygone de Newton}. Si les pentes
de~$(M,F)$ sont les rationnels $\alpha_1\leq\dots\leq\alpha_n$,
compt\'ees avec multiplicit\'es (donc $\dim M=n$), c'est par d\'efinition
l'unique fonction $\Nwt_M:[0;n]\ra\R$ qui est affine
par morceaux, continue, v\'erifie $\Nwt_M(0)=0$ et est de pente
$\alpha_j$ sur l'intervalle $[j;j+1]$.

% \section{Pentes de la cohomologie cristalline et congruences
% de la fonction z\^eta}

Lorsque le corps~$k$ est fini, on a une autre caract\'erisation
des pentes. Supposons pour simplifier que $\sigma^a$ soit l'identit\'e, 
ce qui est v\'erifi\'e dans le cas important o\`u $k=\F_{p^a}$ et $\mathcal V=W(k)$.
L'application $F^a$ est alors $K$-lin\'eaire et la d\'ecomposition 
de Jordan fournit une autre caract\'erisation des pentes :
\begin{prop}
Soit $M$ un F-isocristal et supposons que $\sigma^a=\id$.
Soient $\lambda_1,\dots,\lambda_n$ les valeurs propres de l'application
$K$-lin\'eaire $F^a$. Les pentes de $M$ sont les $\ord_p(\lambda_i)/a$.
\end{prop}

Combinons cette proposition avec la formule des traces de Lefschetz 
en cohomologie rigide :
En effet, si $X$ est un sch\'ema s\'epar\'e de type fini sur
un corps fini $\F_q$, avec $q=p^a$, on a pour tout entier~$n\geq 0$,
\[ \card{ X(\F_q)} = \sum_{i=0}^{2\dim X} (-1)^i \Tr(F^{an}|H^i_\rigc(X/K)). \]
Ainsi, on constate que \emph{minorer} les pentes de la 
cohomologie rigide \`a support propre fournit 
des congruences $p$-adiques pour sa fonction z\^eta.
Pr\'ecis\'ement: contr\^oler la partie de pente~$0$ implique
des congruences modulo~$p$ pour $\card{ X(\F_q)}$,
contr\^oler la partie de pentes~$<1$ des congruences modulo~$q$.

Dans cette veine, il faut citer un r\'esultat fondamental,
d\^u \`a Mazur~\cite{mazur72,mazur73} moyennant des hypoth\`eses
restrictives, et Ogus~\cite[chap.~8]{berthelot-o78} en g\'en\'eral.
Soit $X$ un $k$-sch\'ema propre et lisse ;
sa cohomologie rigide (cristalline en fait) fournit un 
F-isocristal $H^m_\cris(X/W(k))\otimes \Frac W(k)$ dont
nous noterons $\Nwt_X^{(m)}$ le polygone de Newton.
Par ailleurs, $X$ a des nombres de Hodge $h_i=\dim_k
H^{m-i}(X,\Omega^i_{X/k})$.
D\'efinissons le $m$\ieme polygone de Hodge de~$X$, $\Hdg_X^{(m)}$,  
comme la fonction
continue affine par morceaux qui vaut $0$ en~$0$, est de
pente~$0$ sur l'intervalle $[0;h_0]$, de pente~$1$ sur l'intervalle
$[h_0;h_0+h_1]$, etc.

\begin{theo}
On a l'in\'egalit\'e $\Nwt_X^{(m)}\geq \Hdg_X^{(m)}$.
\end{theo}

On en d\'eduit des th\'eor\`emes de type Chevalley-Warning et notamment
une autre d\'emonstration du th\'eor\`eme~\ref{theo.ak} de Ax-Katz
\emph{dans le cas lisse et homog\`ene}. 

% C'est d'ailleurs ce th\'eor\`eme 
% qui avait sugg\'er\'e \`a Katz\footnote{V\'erifier dans~\cite{katz1971}.}
% l'\'enonc\'e pr\'ec\'edent.

\begin{coro}[\`a la Chevalley--Warning]
Soit $X$ une intersection compl\`ete lisse de dimension $d$ dans $\P^{n-1}$,
d\'efinie sur le corps fini~$\F_q$. 
Alors, le nombre de points de $X(\F_{q^s})$ est
\'egal au nombre de points de $\P^d(\F_{q^s})$ modulo $q^{bs}$ o\`u
$b$ est le plus petit entier~$i\geq 0$ tel que $\dim H^{d-i}(\Omega^i) \neq 0$.
{\upshape (Si $d_1,\dots,d_r$ sont les degr\'es des hypersurfaces 
qui d\'efinissent~$X$,
l'entier~$b$ est \'egal au plus petit entier sup\'erieur
ou \'egal \`a \mbox{$(n-\sum d_i)/\max(d_i)$.})}
\end{coro}

\section{Pentes de la cohomologie rigide}

Je donne dans ce paragraphe deux r\'esultats g\'en\'eraux concernant
les pentes de la cohomologie rigide \`a support propre.
Le premier d\'ecrit la partie de pente~$0$, le second pr\'ecise
les pentes possibles.

\begin{theo}\label{lemm.illusie}
Pour tout sch\'ema~$X$, s\'epar\'e et de type fini sur un corps alg\'ebriquement
clos de caract\'eristique~$p>0$, on a un isomorphisme canonique
\[ H^i_\et(X,\Q_p) \otimes K  \simeq H^i_\rigc(X/K)^{0}, \]
o\`u l'exposant~$0$ signifie qu'on consid\`ere la partie de pente~$0$
dans le $F$-isocristal donn\'e par la cohomologie rigide.
\end{theo}
Lorsque $X$ est propre et lisse, ce th\'eor\`eme est d\^u \`a Bloch
(lorsque~$p$ est assez grand) et Illusie (pour tout~$p$, \cite{illusie79},
II, 5.4). Dans le cas g\'en\'eral, c'est un r\'esultat d'\'Etesse et Le Stum
(\cite{etesse-ls1997}, prop.~6.3).

Ils se d\'eduit dans ce cas des propri\'et\'es du complexe de de Rham-Witt via
une g\'en\'eralisation de la suite exacte d'Artin-Schreier
\[ 0 \ra (\Z/p\Z)_X \ra \mathcal O_X \xrightarrow{1-F} \mathcal O_X \ra 0
\]
(suite exacte de faisceaux \'etales sur~$X$).
Plus g\'en\'eralement, si $W_n\mathcal O_X$ d\'esigne le faisceau des vecteurs
de Witt de longueur~$n$ sur~$X$, on a une suite exacte 
de faisceaux \'etales sur~$X$ :
\[ 0 \ra (\Z/p^n\Z)_X \ra W_n\mathcal O_X \xrightarrow{1-F} W_n\mathcal O_X\ra
0 \]
qui induit par passage \`a la cohomologie, puis passage \`a la limite
une suite exacte en cohomologie \'etale:
\[ 0 \ra H^i_\et(X,\Z_p) \ra H^i(X,W\mathcal O_X)
\xrightarrow{1-F} H^i(X,W\mathcal O_X) \ra 0 \]
qui identifie $H^i_\et(X,\Q_p)\otimes K$ au plus grand sous-F-isocristal
de pente~$0$ dans le F-isocristal
$H^i(X,W\mathcal O_X)\otimes_{W(k)} K$.
(Compte tenu du fait que $H^i(X,W\mathcal O_X)$ est pour tout
entier~$i$ un $W(k)$-module de type fini,
la surjectivit\'e de~$1-F$ est une propri\'et\'e g\'en\'erale des~$F$-cristaux,
\cf par exemple~\cite{illusie79}, II, 5.3.)
La th\'eorie du complexe de de Rham-Witt de Bloch et Illusie, et en particulier
la d\'eg\'en\'erescence de la suite spectrale des pentes
(\cite{illusie79}, II, 3.5), implique
que ce dernier espace vectoriel est le plus grand sous-F-isocristal 
de pentes~$[0;1\mathclose[$ dans $H^i_\rig(X/K)$.
Le r\'esultat s'ensuit si $X$ est propre et lisse.

Dans le cas g\'en\'eral, \'Etesse et Le Stum combinent
des suites exactes d'Artin-Schreier, le calcul syntomique
de la cohomologie cristalline (initi\'e par Fontaine et Messing
dans~\cite{fontaine-m87}), la cohomologie cristalline
{\og de niveau variable\fg} et un th\'eor\`eme de Berthelot
selon lequel cette derni\`ere permet de calculer la cohomologie rigide.

\medskip

\begin{theo}
Soit $k$ un corps parfait de caract\'eristique~$p>0$,
$X$ un $k$-sch\'ema de type fini de dimension~$d$.

\begin{enumerate}
\item pour tout~$i$, $H^i_\rigc(X/K)$ est un F-isocristal dont
les pentes appartiennent \`a l'intervalle $[\max(0,i-d),\min(i,d)]$ ;
\item si $X$ est lisse, c'est encore vrai de $H^i_\rig(X/K)$.
\end{enumerate}
\end{theo}
(Dans le cas non lisse, la cohomologie rigide fournit encore
des $F$-isocristaux, mais je ne sais pas ce qu'on peut dire des pentes.)

Le b) est cons\'equence du a) compte tenu de la dualit\'e de Poincar\'e
en cohomologie rigide : il existe un morphisme trace
\[ H^{2d}_\rigc(X/K) \xrightarrow{\Tr} K \]
tel que $\Tr\circ F = p^d F$ qui induit par \emph{cup-produit}
des isomorphismes ($X$ est lisse)
\[ H^i_\rig(X/K) \simeq H^{2d-i}_\rigc (X/K)^\vee (-d) \]
o\`u le \emph{twist} $(-d)$ signifie que le Frobenius est multipli\'e
par $p^d$.
Notons au passage que cela implique l'injectivit\'e des Frobenius
sur la cohomologie rigide (\resp la cohomologie rigide \`a support propre).

Avant de montrer le a), faisons quelques remarques.

1) On peut supposer que le corps $k$ est alg\'ebriquement clos
et que $X$ est r\'eduit.
On peut aussi supposer que $X$ est connexe car la cohomologie
rigide d'une r\'eunion disjointe est la somme directe des cohomologies rigides.

2) Si $X$ est projective et lisse, le th\'eor\`eme de comparaison entre
cohomologies rigide et cristalline implique que $H^i_\rig(X/K)$ 
est un F-isocristal \`a pentes positives ou nulles. 
Par dualit\'e de Poincar\'e, elles sont donc $\leq d$. 
Le th\'eor\`eme de Lefschetz faible
en cohomologie cristalline 
implique alors que les pentes appartiennent \`a l'intervalle
$[0;i]$, donc \`a l'intervalle $[i-d;i]$ via la dualit\'e de Poincar\'e.

3) Si $U$ est un ouvert dense de~$X$ et $Z$ le ferm\'e compl\'ementaire,
la suite exacte longue d'excision 
\[  \ra H^i_\rigc(Z/K) \ra H^i_\rigc(U/K) \ra H^{i}_\rigc(X/K) \ra \]
et l'hypoth\`ese que l'assertion a) est v\'erifi\'ee en dimension $<\dim X$
montrent qu'il est \'equivalent de d\'emontrer l'\'enonc\'e a) pour $X$ et pour $U$.

4) Si $U'\ra U$ est un rev\^etement \'etale de $k$-vari\'et\'es lisses,
la cohomologie rigide \`a support propre de $U'$ admet celle de~$U$
comme facteur direct, si bien qu'il suffit de d\'emontrer l'assertion a)
pour $U'$.

D\'emontrons maintenant a) par r\'ecurrence sur la dimension de~$X$.
D'apr\`es le th\'eor\`eme d'alt\'eration de de Jong (\cite{dejong96}, th.~4.1),
il existe un ouvert dense $U\subset X$, une $k$-vari\'et\'e projective et lisse
$X'$, un ouvert $U'\subset X'$ et un rev\^etement \'etale $U'\ra U$.

D'apr\`es 2), l'assertion a) est vraie pour $X'$. D'apr\`es 3),
elle est donc vraie pour $U'$ et la remarque 4) entra\^{\i}ne sa v\'eracit\'e
pour~$U$, donc aussi pour~$X$ gr\^ace \`a 3).

\begin{rema}[R\'ef\'erences bibliographiques]
La d\'emonstration est celle sugg\'er\'ee par Berthelot
dans~\cite{berthelot1997}, remarque~3.9 et se trouve aussi
dans l'article~\cite{chiarellotto-ls1999} de B.~Chiarellotto et B.~Le Stum.
En suivant cette approche, ces auteurs ont aussi
\'elucid\'e la structure des poids (c'est-\`a-dire 
des valeurs absolues archim\'ediennes
des valeurs propres de Frobenius) sur la cohomologie rigide \`a support
propre d'une vari\'et\'e sur un corps fini
(cf.~ \cite{chiarellotto1998} et \cite{chiarellotto-ls1999b}).
Ils doivent faire usage du th\'eor\`eme de Katz-Messing :
dans ~\cite{katz-m74}, ces derniers
d\'eduisent des conjectures de Weil 
et du th\'eor\`eme de Lefschetz difficile
en cohomologie $\ell$-adique les th\'eor\`emes correspondants
en cohomologie cristalline. Signalons aussi que Kedlaya a r\'ecemment
adapt\'e \`a la cohomologie rigide (voir~\cite{kedlaya2002b}).
la d\'emonstration par Laumon des conjectures de Weil.

On peut en fait d\'emontrer un analogue de ce th\'eor\`eme sur un corps
fini, via la cohomologie \'etale $\ell$-adique, et c'est ainsi
que proc\`ede Ekedahl dans~\cite{ekedahl1983}.
Cela n\'ecessite de montrer au pr\'ealable que les valeurs propres de Frobenius
sur la cohomologie $\ell$-adique \`a support propre
sont des entiers alg\'ebriques, ce qui est fait par Deligne 
(\S\,5 de l'expos\'e~\cite{katz1973b}). On peut alors, par un d\'evissage
analogue, \'etudier leurs valuations $p$-adiques. 
 
Enfin, signalons un article de M.~Kim~\cite{kim2002} dans lequel
la cohomologie rigide est remplac\'ee par celle d'un complexe
de de Rham-Witt \`a p\^oles logarithmiques.
\end{rema}

\section{Les th\'eor\`emes d'Ekedahl et Esnault}

Dans ce paragraphe, je pr\'esente deux th\'eor\`emes dus \`a
T.~Ekedahl~\cite{ekedahl1983}
et H.~Esnault~\cite{esnault2003} qui permettent
de contr\^oler la partie de pentes $<1$ 
dans la cohomologie rigide.
Leur utilit\'e appara\^{\i}tra aux paragraphes suivants,
lorsque nous d\'eduirons de cette partie de pente $<1$ des renseignements
g\'eom\'etriques.
Dans~\cite{kahn2002}, B.~Kahn red\'emontre ces \'enonc\'es
analogues \`a l'aide des \emph{motifs birationnels},
notion qu'il a introduite avec R.~Sujatha. 

Si $(M,F)$ est un F-isocristal et si $\alpha\in\Q$,
notons $M^{>\alpha}$ et $M^{<\alpha}$ les plus grands
sous-F-isocristaux de $M$ dont les pentes soient respectivement
strictement sup\'erieures et inf\'erieures \`a $\alpha$.

\begin{theo}[Ekedahl] 
\label{theo.ekedahl}
Soit $k$ un corps parfait de caract\'eristique positive.
Soit $X$ et $Y$ deux $k$-sch\'emas de type fini, int\`egres
de m\^eme dimension~$d$.
\begin{enumerate}\def\theenumi{\alph{enumi}}
\item S'il existe une application rationnelle dominante $X\ra Y$,
il existe pour tout $i$, une injection de F-isocristaux (canonique) 
de $H^i_{\rigc}(Y)^{>d-1}$  dans $H^i_{\rigc}(X)^{>d-1}$ ;
\addtocounter{enumi}{-1}\def\theenumi{\alph{enumi}$'$}
\item S'il existe une application rationnelle dominante $X\ra Y$
et si $X$ et $Y$ sont lisses, $H^i_\rig(Y)^{<1}$ s'injecte
canoniquement dans $H^i_\rig(X)^{<1}$.
\def\theenumi{\alph{enumi}}
\item Si $X$ et $Y$ sont birationnels, alors pour tout~$i$,
$H^i_\rigc(X)^{>d-1}$ et $H^i_{\rigc}(X)^{>d-1}$ sont des F-isocristaux
isomorphes.
\addtocounter{enumi}{-1}\def\theenumi{\alph{enumi}$'$}
\item Si $X$ et $Y$ sont lisses et birationnelles, alors
pour tout~$i$, les F-isocristaux $H^i_\rig(X)^{<1}$ et $H^i_\rig(Y)^{<1}$
sont isomorphes.
\end{enumerate}
\end{theo}

Les \'enonc\'es a$'$) et b$'$) se d\'eduisent des \'enonc\'es a) et b)
par dualit\'e de Poincar\'e. Il est par ailleurs clair que a) implique b),
ce qui nous ram\`ene \`a d\'emontrer l'\'enonc\'e a).

Il existe des ouverts lisses et denses $U\subset X$ et $V\subset Y$ 
sur lesquelles l'application rationnelle $X\ra Y$ induit
un morphisme fini et plat $f\colon U\ra V$. Les suites exactes d'excision
associ\'ees aux ouverts~$U$ et $V$ et l'\'etude des pentes de la cohomologie
rigide \`a support propre impliquent que l'on a des isomorphismes
(induits par les inclusions) 
\[ H^i_\rigc(U/K)^{>d-1} \ra H^i_\rigc(X/K)^{>d-1}, \qquad
   H^i_\rigc(V/K)^{>d-1} \ra H^i_\rigc(Y/K)^{>d-1}. \]
Par ailleurs, \`a $f$ est attach\'e deux morphismes de F-isocristaux
\[ f_*\colon H^i_\rigc(U/K)\ra H^i_\rigc(V/K) \quad\text{et}\quad
  f^*\colon H^i_\rigc(V/K) \ra H^i_\rigc(U/K) \]
tels que $f_*\circ f^*=\deg(f)$. Il en r\'esulte que $f^*$ est injectif 
et le r\'esultat s'en d\'eduit.

\begin{theo}[Esnault]
\label{theo.esnault}
Soit $k$ un corps de caract\'eristique $p>0$.
Soit $X$ une vari\'et\'e propre et lisse telle que, notant
$\Omega$ une cl\^oture alg\'ebrique du corps des fonctions $k(X)$,
$\CH_0(X_\Omega)_\Q=\Q$.
Alors, pour tout $i>0$, $H^i_\rig(X/K)^{<1}=0$.
\end{theo}
On peut supposer que $k$ est alg\'ebriquement clos.
D'apr\`es un th\'eor\`eme de Bloch
(\cite{bloch1980,bloch1980b}, 
je rappelle aussi l'argument plus bas),
il existe un ouvert dense~$U$ de~$X$ et un point $x_0\in X(k)$
tel que sur~$U\times X$,   notant $\Delta_X$ la diagonale de $X\times X$,
\begin{equation}\label{eq.dec}
 \Delta_X = U\times [x_0] 
\quad
\text{dans $\CH_d(U\times X)_\Q$.} 
\end{equation}
Autrement dit,  le graphe $\Gamma_j\subset U\times X$
de l'immersion ferm\'ee $j\colon U\ra X$ est lin\'eairement
\'equivalent \`a $U\times [x_0]$ dans $\CH_d(U\times X)_\Q$.

Remarquons qu'une classe de cohomologie
$\alpha\in H^{2d}_\rig(U\times X)$ d\'efinit une correspondance
sur la cohomologie rigide \`a support propre :
\[ \alpha_*\colon H^i_\rigc(U) \xrightarrow{q^*} H^i_\rigc(U\times X)
 \xrightarrow{\cap\alpha} H^{i+2d}_\rigc(U\times X) \xrightarrow{p_*}
H^i_\rigc(X). \]
L'application $q^*$ existe car $X$ est propre, l'accouplement
$\cap\alpha$ vient de celui entrep cohomologie rigide \`a support
propre et cohomologie rigide (d\'efini dans~\cite{berthelot1997b})
et l'application $p_*$  est la transpos\'ee, via la dualit\'e de Poincar\'e,
de l'application $p^*$ sur la cohomologie rigide.
Le graphe de $j$ a une classe $\gamma(\Gamma_j)$
dans $H^{2d}_\rig(U\times X)$, construite dans~\cite{petrequin2003}
et l'on a pour tout $i$,
\[ \gamma(\Gamma_j)_* = j_* \colon H^i_\rigc(U)\ra H^i_\rigc(X). \]
Mais la d\'ecomposition~\eqref{eq.dec} affirme que modulo
l'\'equivalence rationnelle, $\Gamma_j=p^*[x_0]$,
donc, en passant aux classes de cycles (\cite{petrequin2003}, prop.~6.10),
\[ \gamma(\Gamma_j) = p^* \gamma([x_0]) , \qquad \gamma([x_0])
        \in H^{2d}_\rig(X) \]
Par suite, pour toute classe $\xi\in H^i_\rigc(U)$,
\begin{align*}
 j_*(\xi) &= \gamma(\Gamma_j)_*(\xi) = p_*(q^*(\xi)\cap p^*(\gamma([x_0])) \\
& =  p_*(q^*(\xi)) \cap \gamma([x_0]).
\end{align*}
Comme $p_*(q^*(\xi))\in H^{i-2d}_\rigc(X)$, il est nul pour $i<2d$
et finalement, les homomorphismes $j_*\colon H^i\rigc(U)\ra H^i_\rigc(X)$
sont nuls pour $0\leq i<2d$. (En degr\'e $2d$, c'est un isomorphisme.)
D'apr\`es le th\'eor\`eme~\ref{theo.ekedahl},
ils induisent des isomorphismes  sur les parties de pentes $>d-1$.
Par cons\'equent, pour tout $i$, le F-isocristal $H^i_\rigc(X)^{>d-1}$
est nul.

Par dualit\'e de Poincar\'e, $H^i_\rig(X)^{<1}=0$ si $i>0$.

\medskip

Donnons pour terminer la d\'emonstration de l'existence d'une
d\'ecomposition~\eqref{eq.dec}.
L'hypoth\`ese est que $\CH_0(X_\Omega)_\Q=\Q$, o\`u $\Omega$
est une cl\^oture alg\'ebrique de $k(X)$. On suppose
toujours que le corps~$k$ est alg\'ebriquement clos.
Soit $x_0$ un point g\'eom\'etrique
de $X$ et soit $\eta$ son point g\'en\'erique; ils d\'efinissent deux points
de $X(k(X))$ et leur diff\'erence~$\alpha$ est un $0$-cycle sur $X_{k(X)}$.
Par hypoth\`ese, leur diff\'erence est de torsion sur $X_\Omega$,
donc sur une extension finie~$K$ de $k(X)$. En prenant des normes,
on voit que $\alpha$ est d\'ej\`a de torsion dans $\CH_0(X_{k(X)})$.
Cela signifie qu'il existe un ouvert~$U$ de~$X$ tel que le $d$-cycle
$\bar\alpha=X\times [x_0] - \Delta_X$ est de torsion sur $U\times X$.

\bigskip

Rappelons maintenant quelques conditions qui affirment
que le groupe $\CH_0(X_\Omega)_\Q$ est \'egal \`a $\Q$.

Une classe importante de vari\'et\'es o\`u elle est v\'erifi\'ee
est fournie par les vari\'et\'es (g\'eom\'etriquement) rationnellement
connexes : 
il s'agit de vari\'et\'es propres~$X$ pour lesquelles deux points g\'en\'eraux
(disons sur une extension du corps de base de cardinal non d\'enombrable)
sont joints par une courbe rationnelle, c'est-\`a-dire par l'image
d'un morphisme $\P^1\ra X$.
Les vari\'et\'es faiblement unirationnelles
sont \'evidemment rationnellement connexes ; l'\'enonc\'e b)
redonne ainsi l'un des r\'esultats principaux de la note~\cite{ekedahl1983}.

Plus g\'en\'eralement, une vari\'et\'e~$X$ est dite \emph{rationnellement
connexe par cha\^{\i}nes} si deux points quelconques 
sont joints par une cha\^{\i}ne de courbes rationnelles.
D'apr\`es un th\'eor\`eme d\^u ind\'ependamment \`a F.~Campana
et J.~Koll\'ar, Y.~Miyaoka et S.~Mori,
les vari\'et\'es de Fano, c'est-\`a-dire les vari\'et\'es projectives
lisses dont le fibr\'e anticanonique est ample, sont
rationnellement connexes par cha\^{\i}nes.
(Voir~\cite{debarre2001}, prop.~5.16.)
En caract\'eristique z\'ero, une vari\'et\'e propre et lisse
qui est rationnellement connexe par cha\^{\i}nes est en fait
rationnellement connexe. 

On peut d\'emontrer que si $X$ est rationnellement connexe par
cha\^{\i}nes, deux points \emph{quelconques} 
de~$X$ sont reli\'es par une cha\^{\i}ne de courbes
rationnelles. Ainsi, \mbox{$\CH_0(X_\Omega)_\Q=\Q$.}

\begin{rema}
Si $\dim X\leq 3$, et si $X$ est s\'eparablement unirationnelle,
ou bien de Fano, on sait que
les groupes de cohomologie $H^i(X,\mathcal O_X)$ sont nuls
pour $i\geq 1$
(Nygaard~\cite{nygaard1978}, Shepherd-Barron~\cite{shepherd-barron1997}).
On peut alors en d\'eduire que pour $i\geq 1$,
$H^i(X,W\mathcal O_X)=0$, d'o\`u une autre approche au
th\'eor\`eme~\ref{theo.esnault}.
\end{rema}

\begin{rema}
Dans~\cite{esnault2003}, H.~Esnault \'enonce son th\'eor\`eme
sous l'hypoth\`ese, apparemment plus forte,
que $\CH^0(X_\Omega)=\Z$ ; elle est en fait \'equivalente.
Notons $\CH^0(X_\Omega)^0$  le noyau de l'application
degr\'e $\CH^0(X_\Omega)\ra \Z$ et  soit $\alb\colon X\ra A$ l'application
d'Albanese de~$X_\Omega$.
Comme $A$ est engendr\'ee par l'image de~$X$, 
l'application
$\CH^0(X_\Omega)^0 \ra A(\Omega)$ d\'eduite de $\alb$
est surjective. Si $\CH^0(X_\Omega)_\Q=\Q$, $\CH^0(X_\Omega)^0$
est un groupe de torsion, $A(\Omega)$ aussi, et donc $A=0$.
Par ailleurs, un th\'eor\`eme de Ro\u\i tman~\cite{roitman1980} compl\'et\'e par
Milne~\cite{milne1982}
affirme que l'application $\CH^0(X_\Omega)^0\ra A(\Omega)$ induit
un isomorphisme sur les sous-groupes de torsion. Il en r\'esulte
que $\CH^0(X_\Omega)^0$ est sans torsion, donc nul.
Ainsi, $\CH^0(X_\Omega)=\Z$.
\end{rema}

% \part{Fromage et dessert}

\section{Retour sur la fonction z\^eta}

Comme je l'ai d\'ej\`a mentionn\'e au \S\,5, la formule des traces
de Lefschetz en cohomologie rigide \`a support propre 
montre que la partie de
la cohomologie de pente $<\alpha$ fournit une congruence
modulo $q^\alpha$ pour le nombre de points rationnels.

Par suite, les th\'eor\`emes~\ref{theo.ekedahl} et~\ref{theo.esnault}
impliquent le th\'eor\`eme suivant:
\begin{theo}
Soit $X$ et $Y$ deux vari\'et\'es propres, lisses et g\'eom\'etriquement
connexes sur le corps fini $\F_q$. Soit $\Omega$
une cl\^oture alg\'ebrique de $\F_q(X)$.
\begin{enumerate}\def\theenumi{\alph{enumi}}
\item Si $X_\Omega$ et $Y_\Omega$ sont birationnelles, alors 
$\card{ X(\F_q)}\equiv \card{ Y(\F_q)} \pmod q$ ;
\item Si $\CH_0(X_\Omega)_\Q=\Q$, alors $\card{ X(\F_q)}\equiv 1\pmod q$.
\end{enumerate}
\end{theo}

M\^eme si l'\'enonc\'e a) ne figure pas explicitement dans la
note~\cite{ekedahl1983}, c'en est une cons\'equence imm\'ediate.
D'un autre c\^ot\'e, il est \'evident si $X$ et $Y$ sont li\'ees 
par un \'eclatement de centre lisse.
Il le serait plus g\'en\'eralement pour tout couple de vari\'et\'es birationnelles,
si l'on disposait d'un th\'eor\`eme de factorisation faible
en caract\'eristique positive. 
Dans leur article r\'ecent~\cite{lachaud-p2000}, G.~Lachaud et M.~Perret
ont fait marcher cette approche en dimension~3.

Puisque les vari\'et\'es de Fano sur un corps
alg\'ebriquement clos
sont rationnellement connexes par cha\^{\i}nes,
il en r\'esulte ainsi le th\'eor\`eme, 
conjectur\'e par Lang et Manin~\cite{manin1993}:
\begin{coro}[Esnault]
Soit $X$ une vari\'et\'e de Fano sur un corps fini $\F_q$.
Alors, $\card{ X(\F_q)}\equiv 1\pmod q$. Il est en particulier
non nul.
\end{coro}

On en d\'eduit aussi le r\'esultat:
\begin{prop}
Soit $X$ une vari\'et\'e propre, lisse et g\'eom\'etriquement connexe sur
un corps fini  $\F_q$. Supposons que $X$ soit g\'eom\'etriquement
domin\'ee par une $k$-vari\'et\'e $Y$, propre, lisse et connexe,
de m\^eme dimension,
telle que $H^i_\et(Y, \Q_p)=0$ si $i\neq 0$.
Alors, pour toute extension finie $\F$ de~$\F_q$,
$\card{ X(\F)}\equiv 1\pmod p$.
\end{prop}

\begin{rema}
Le th\'eor\`eme d'Esnault peut \^etre mis en parall\`ele avec
plusieurs r\'esultats r\'ecents.
Soit $k$ le corps des fonctions d'une courbe projective
et lisse sur un corps alg\'ebriquement clos.
Graber, Harris, Starr~\cite{graber-h-s2003}, 
et de Jong et Starr~\cite{dejong-s2003}
ont montr\'e qu'une $k$-vari\'et\'e projective, lisse qui  est
s\'eparablement rationnellement connexe a un point rationnel.
Un tel corps~$k$, de m\^eme qu'un corps fini, est $\mathrm C_1$,
donc de dimension cohomologique au plus~$1$.
Cependant, Colliot-Th\'el\`ene et Madore ont construit
dans~\cite{colliot-thelene-m2003} une surface cubique
sur un corps $p$-adique qui n'a pas de point rationnel.
La question de savoir si une vari\'et\'e s\'eparablement rationnellement
connexe (voire rationnellement connexe par cha\^{\i}nes)
sur un corps $\mathrm C_1$ admet un point rationnel
reste ouverte.
\end{rema}

\begin{rema}
Contrairement au th\'eor\`emes d'Ax et Katz, la d\'emonstration 
du th\'eor\`eme~\ref{theo.esnault} n\'ecessite une hypoth\`ese de lissit\'e
Bloch, Esnault et Levine~\cite{bloch-e-l2003}
ont propos\'e de remplacer  la d\'ecomposition de la diagonale
dans le groupe de Chow (formule~\eqref{eq.dec}) par
une d\'ecomposition analogue dans un groupe de \emph{cohomologie motivique}
ad\'equat. Leur condition entra\^{\i}ne une minoration des pentes de la cohomologie 
rigide et ils ont montr\'e qu'elle est v\'erifi\'ee
dans le cas des hypersurfaces de degr\'e~$d\leq n$
\'eventuellement singuli\`eres de l'espace projectif~$\mathbf P^n$.
\end{rema}

 Les vari\'et\'es singuli\`eres ne sont pas couvertes par le
th\'eor\`eme~\ref{theo.esnault}, alors que
\section{Simple connexit\'e de certaines vari\'et\'es}

On peut aussi appliquer ces consid\'erations pour \'etudier
le groupe fondamental de certaines vari\'et\'es alg\'ebriques,
notamment les vari\'et\'es unirationnelles.

Rappelons que J-P.~Serre a d\'emontr\'e dans~ \cite{serre1959}
que le groupe fondamental d'une telle vari\'et\'e est trivial,
pourvu que le corps de base soit \emph{de caract\'eristique z\'ero.}

\begin{lemm}
Soit $X$ une vari\'et\'e projective lisse, g\'eom\'etriquement connexe,
sur un corps de caract\'eristique z\'ero.
On suppose que $X$ est unirationnelle, ou que $X$ est de Fano,
ou, $\Omega$ d\'esignant la cl\^oture alg\'ebrique du corps
des fonctions de~$X$, que $\CH^0(X_\Omega)_\Q=\Q$.
Alors, $\chi(X,\mathcal O_X)=1$.
\end{lemm}
Dans le premier cas, on a en effet $H^0(X,\Omega^i_X)=0$, pour $i>0$,
comme on le voit
en tirant une $i$-forme de $X$ \`a l'espace projectif  : elle y sera r\'eguli\`ere
hors d'un lieu de codimension~2, donc partout, donc nulle.
 En caract\'eristique
z\'ero, cet espace a m\^eme dimension que $H^i(X,\mathcal O_X)$, d'o\`u
l'assertion.
Dans le second cas, si $i>0$, $H^i(X,\mathcal O_X)$ est dual de
$H^{d-i}(X,\omega_X^{-1})$, donc est nul par le th\'eor\`eme d'annulation
de Kodaira.

Le dernier cas se d\'emontre en  d\'ecomposant la diagonale mais
en th\'eorie de Hodge. Par le m\^eme argument que dans la preuve
du th\'eor\`eme~\ref{theo.esnault}, il existe un ouvert dense~$U$ de~$X$
tel que l'application $H^i_{DR}(X)\ra H^i_{DR}(U)$ soit nulle pour $i>0$.
D'autre part, la th\'eorie de Hodge mixte de P.~Deligne
fournit une factorisation
\[ 
      H^i_{DR}(X) \rightarrow H^i_{DR}(U) \rightarrow H^i(X,\mathcal O_X),
\]
l'application compos\'ee \'etant surjective
(\cf\cite{deligne1972}, (3.2.13), (\emph{ii})).
Il en r\'esulte que $H^i(X,\mathcal O_X)=0$ pour $i>0$.

\begin{coro}\label{coro.car0-triv}
\emph{(En caract\'eristique z\'ero.)}
Une vari\'et\'e propre et lisse qui est
unirationnelle, ou de Fano, ou rationnellement connexe par cha\^{\i}nes,
n'a pas de rev\^etement \'etale fini non trivial.
\end{coro}
Si $f\colon Y\ra X$ est un tel rev\^etement, remarquons  que l'hypoth\`ese
implique que $Y$ est aussi unirationnelle (\resp de Fano).
Par suite, on a $\chi(Y,\mathcal O_Y)=1$.
Or, le th\'eor\`eme de Riemann-Roch implique que $f_*\ch(\mathcal O_Y)
=\deg(f)\ch(\mathcal O_X)$, si bien que $\chi(Y,\mathcal
O_Y)=\deg(f)\chi(X,\mathcal O_X)$.
N\'ecessairement, $\deg(f)=1$.

\begin{rema}
a) Si le corps de base est~$\C$, le groupe fondamental topologique
de telles vari\'et\'es est aussi trivial.

b) 
Soit $X$ une surface d'Enriques sur le corps des nombres
complexes. Bloch, Kas et Lieberman montrent 
dans~\cite{bloch-k-l1976} que sur~$X$, tout z\'ero-cycle 
et de degr\'e nul est rationnellement \'equivalent \`a~$0$.
L'hypoth\`ese $\CH_0(X_\Omega)_\Q=\Q$
ne suffit donc pas \`a assurer la validit\'e du corollaire pr\'ec\'edent.
\end{rema}

\medskip

En caract\'eristique~$p>0$, les choses sont plus compliqu\'ees.
Tout d'abord, il existe des surfaces unirationnelles non
simplement connexes (Shioda, \cite{shioda1974}) : si $p\neq 5$
et $p\not\equiv 1\pmod 5$, la surface d'\'equation
\[ X_0^5+X_1^5+X_2^5+X_3^5 = 0 \]
dans $\P^3$ est unirationnelle mais poss\`ede une action libre
du groupe des racines 5\iemes de l'unit\'e (par $X_i\ra \zeta^i X_i$).
La surface
de Godeaux obtenue par quotient est alors
unirationnelle mais n'est pas simplement connexe : son groupe
fondamental est $\Z/5\Z$.
Cet exemple montre aussi qu'en caract\'eristique positive,
une vari\'et\'e rationnellement connexe par cha\^{\i}nes~$X$ ne
v\'erifie pas forc\'ement $H^i(X,\mathcal O_X)=0$.

\medskip

Nous allons d\'emontrer une g\'en\'eralisation d'un r\'esultat d'Ekedahl
(\cite{ekedahl1983})
qui repose sur une formule d'Euler-Poincar\'e en cohomologie \'etale
$p$-adique, d\'emontr\'ee par R.~Crew dans~\cite{crew1984}.
\begin{prop}\label{prop.crew}
Soit $k$ un corps alg\'ebriquement clos
de caract\'eristique~$p$
et soit $f\colon Y\ra X$ un rev\^etement \'etale galoisien
de $k$-sch\'emas s\'epar\'es de type fini, de degr\'e une puissance de~$p$.
On a alors les formules suivantes entre caract\'eristiques 
d'Euler-Poincar\'e en cohomologie \'etale \`a support propre :
\[ \chi_\etc(Y,\Q_p) = \deg(f) \chi_\etc(X,\Q_p)
\quad\text{et}\quad
     \chi_\etc(Y,\Z/p\Z) = \deg(f) \chi_\etc(X,\Z/p\Z).  \]
\end{prop}
En consid\'erant la suite exacte
\[ 0 \ra \Z_p\ra \Z_p \ra \Z/p \ra 0 \]
et en utilisant le fait que la cohomologie \'etale \`a support
propre \`a coefficients dans $\Z_p$ est repr\'esent\'ee par
un complexe parfait de $\Z_p$-modules,
on d\'emontre que 
\[ \chi_\etc(X,\Q_p)=\chi_\etc(X,\Z/p\Z) 
\quad\text{et}\quad
\chi_\etc(X,\Q_p)=\chi_\etc(X,\Z/p\Z), \]
si bien que la seconde formule implique la premi\`ere.

Le faisceau \'etale $f_*(\Z/p\Z)$ sur $X$ est localement constant
et correspond \`a une repr\'esentation de son groupe fondamental 
sur un $\F_p$-espace vectoriel de dimension~$d=\deg(f)$,
repr\'esentation qui provient d'une repr\'esentation de $\Gal(Y/X)$.
Comme ce groupe est un $p$-groupe, cette repr\'esentation est extension
successive de repr\'esentations triviales. Autrement dit,
$f_*(\Z/p\Z)$ est extension successive de~$d$ faisceaux $\Z/p\Z$.
Par suite
\[ \chi_\etc(Y,\Z/p\Z) = \chi_\etc(X,f_*(\Z/p\Z))=d \chi_\etc(X,\Z/p\Z). \]

\bigskip

D'apr\`es le lemme~\ref{lemm.illusie},
si $X$ est unirationnelle, \resp Fano, \resp rationnellement connexe
par cha\^{\i}nes,
sur un corps alg\'ebriquement clos $k$ de caract\'eristique~$p$, 
on a $\chi_\etc(X,\Q_p)=1$.
Mais un rev\^etement \'etale d'un tel~$X$ est aussi
unirationnel (\resp Fano, \resp rationnellement connexe par cha\^{\i}nes).
L'argument utilis\'e dans la preuve du corollaire~\ref{coro.car0-triv}
permet alors de montrer le th\'eor\`eme suivant, d\^u \`a Ekedahl
dans le cas unirationnel.
\begin{prop}\label{prop.p'}
Soit $k$ un corps alg\'ebriquement clos de caract\'eristique~$p>0$
et soit $X$ une vari\'et\'e propre, lisse sur~$k$.

Si $X$ est Fano, ou si $X$ est rationnellement connexe par cha\^{\i}nes,
son groupe fondamental est un groupe fini d'ordre premier \`a~$p$.
\end{prop}

Avant de d\'emontrer ce r\'esultat, rappelons que 
les vari\'et\'es dites \emph{s\'eparablement rationnellement connexes}
sont simplement connexes.
(Cela r\'esulte du th\'eor\`eme de de Jong et Starr, \cf
l'expos\'e~\cite{debarre2002} d'O.~Debarre).

Par ailleurs, on peut d\'emontrer qu'une vari\'et\'e 
propre, d\'efinie sur un corps alg\'ebriquement clos
qui est normale et rationnellement connexe par cha\^{\i}nes,
a un groupe fondamental fini
(Campana~\cite{campana1995} en caract\'eristique z\'ero.
Voir~\cite{acl2003} pour le cas g\'en\'eral.)

% Si $k$ est alg\'ebriquement clos, on d\'emontre en effet
% (\cite{debarre2001}, p.~89) qu'une $k$-vari\'et\'e propre est
% rationnellement connexe si et seulement s'il existe un point
% $x_0\in X(k)$, un $k$-sch\'ema $M$
% et une application dominante~$\gamma\colon \P^1\times M \ra X$
% telle que pour tout~$m\in M$, $\gamma(0,m)=x_0$.
% Par suite, d'une part l'image de
% $ \pi_1(\P^1\times M)$ dans $\pi_1(X)$ est d'indice fini
% et d'autre part, cette image est nulle comme le montre la factorisation
% \[ \pi_1(\P^1\times M)\simeq \pi_1(\{0\}\times M)  \ra \pi_1(x_0)
% \ra\pi_1(X). \]
% (En caract\'eristique positive, 
% le premier isomorphisme r\'esulte de ce que $\P^1$ est propre.)

D\'emontrons maintenant la proposition~\ref{prop.p'}.
Soit $Y\ra X$ un rev\^etement \'etale galoisien de groupe~$G$.
Si $X$ est de Fano (\resp rationnellement connexe par cha\^{\i}nes),
notons que $Y$ l'est aussi. Soit $P$ un $p$-sous-groupe de Sylow de~$G$,
de sorte que le rev\^etement $Y\ra Y/P$ est \'etale galoisien de groupe~$P$.
Les vari\'et\'es $Y$, $Y/P$ v\'erifient
$\chi_\et(Y,\Q_p)=\chi_\et(Y/P,\Q_p)=1$ et la formule
de Crew affirme que $\chi_\et(Y,\Q_p)=\card{P} \chi_\et(Y/P,\Q_p)$.
On a donc $P=\{1\}$.

\clearpage

\bibliographystyle{smfplain}
\bibliography{aclab,biblio,acl}

\addressindent.47\textwidth
\end{document}